\documentclass[11pt,letterpaper,twoside]{article}
\usepackage[centertags]{amsmath}
\usepackage{graphicx,indentfirst,amsmath,amsfonts,amssymb,amsthm,newlfont}
\font\Goth=yinitas scaled \magstep0

\newcommand{\Gth}[1]{\lower2mm\hbox{\Goth #1}}

\usepackage[latin1]{inputenc}
\def\al{\alpha}

\def\de{\delta}

\def\l1{{\lambda}_1}

\newcommand{\f}{\frac}

\def\x1{{\xi }_{xx}}
\def\x2{{\xi }_{yy}}
\def\x3{{\xi }_{xy}}

\def\e1{{\eta }_{xx}}
\def\e2{{\eta }_{yy}}
\def\e3{{\eta }_{xy}}

\newcommand{\ds}{\displaystyle }
\newtheorem{theorem}{Theorem}
\newtheorem{lemma}{\bf Lemma}
\newtheorem{corollary}{Corollary}
\newcommand{\beqn}{\begin{eqnarray*}}
\newcommand{\eeqn}{\end{eqnarray*}}
\newcommand{\beqnn}{\begin{eqnarray}}
\newcommand{\eeqnn}{\end{eqnarray}}
\newcommand{\p}{\partial}
\newcommand{\bb}{\begin{equation}}
\newcommand{\ee}{\end{equation}}
\newcommand{\ba}{\begin{array}}
\newcommand{\ea}{\end{array}}
\newcommand{\R}{\mathbb{R}}

\newcommand{\lb}{\Delta_{g}}
\newcommand{\s}{\mathbb{S}^{2}\times\R}

\begin{document}
\pagenumbering{arabic}
\title{\huge \bf On the paper ``Symmetry analysis of wave equation on sphere'' by H. Azad and M. T. Mustafa}
\author{\rm \large Igor Leite Freire \\
\\
\it Centro de Matemática, Computação e Cognição\\ \it Universidade Federal do ABC - UFABC\\ \it 
Rua Catequese, $242$,
Jardim,
$09090-400$\\\it Santo André, SP - Brasil\\
\rm E-mail: igor.freire@ufabc.edu.br}
\date{\ }
\maketitle
\vspace{1cm}
\begin{abstract}
Using the scalar curvature of the product manifold $\mathbb{S}^{2}\times\R$ and the complete group classification of nonlinear Poisson equation on (pseudo) Riemannian manifolds, we extend the previous results on symmetry analysis of homogeneous wave equation obtained by H. Azad and M. T. Mustafa [H. Azad and M. T. Mustafa, Symmetry analysis of wave equation on sphere, J. Math. Anal. Appl., 333 (2007) 1180--1188] to nonlinear Klein-Gordon equations on the two-dimensional sphere.\end{abstract}
\vskip 1cm
\begin{center}
{2000 AMS Mathematics Classification numbers:\vspace{0.2cm}\\
76M60, 58J70, 35A30, 70G65\vspace{0.2cm} \\
Key words: Lie point symmetry, Noether symmetry, conservation laws, wave equations on the sphere, scalar curvature}
\end{center}
\pagenumbering{arabic}
\newpage

\section{Introduction}

\

In a previous work, Azad and Mustafa \cite{A} considered the Lie point symmetries of the homogeneous wave equation induced by the $2-$Sphere $\mathbb{S}^{2}$ metric
\bb\label{1.1}
u_{tt}=u_{xx}+(\cot{x})\;u_{x}+\f{1}{\sin^{2}{x}}\;u_{yy}.
\ee

Equation (\ref{1.1}) is a particular case of the general equation
\bb\label{1.2}
\lb u+f(u)=0,
\ee
where 
$$\lb u=\f{1}{\sqrt{g}}\f{\p}{\p x^{i}}(\sqrt{g}g^{ij}\f{\p u}{\p x^{j}})=
g^{ij}\nabla_{i}\nabla_{j}u=\nabla^{j}\nabla_{j}u=\nabla_{i}\nabla^{i}u,$$
where $\lb$ is the Laplace-Beltrami operator on an arbitrary (pseudo) Riemannian manifold $(M^{n},g)$ and
${\nabla }_i$ is the covariant derviative corresponding to the
Levi-Civita connection and the Einstein summation
convention over repeated indices is understood. 

Equation (\ref{1.2}) covers Poisson and Klein-Gordon semilinear equations, depending on if $(M^{n},g)$ is a Riemannian or a pseudo-Riemannian manifold, respectively. Equation (\ref{1.1}) can be obtained from (\ref{1.2}) taking on $\s$ the metric 
\bb\label{met}ds^{2}=dt^{2}-dx^{2}-\sin^{2}{x}\;dy^{2}\ee
and $f(u)=0$.
 
We shall denote the product manifold $\s$ endowed with the metric (\ref{met}) as $(\s,g)$.

Group classification of equations with coefficients depending on metric tensor on especific Riemannian manifolds are well-known. See \cite{A,gc,yi1,yii1,yii2,ib}. 

The Lie point symmetries of equation (\ref{1.2}) on flat manifolds, with some functions $f(u)$ are performed in \cite{ib}. In \cite{gc,ds,cl,i} the Lie point symmetries, the Noether symmetries and the conservation laws of the Kohn-Laplace equations were studied. In \cite{yi1} the symmetry analysis of equation (\ref{1.2}) was carried out on an arbitrary (pseudo) Riemannian manifold. The Lie symmetries of the Poisson equation with Euclidean metric are well known, see \cite{sv}. The group classification of the Poisson equation on the Hyperbolic plane with metric of Klein's model of Lobachevskian geometry and in $\mathbb{S}^{2}$ was carried out in \cite{yii1,yii2} and \cite{yii2}, respectively.

In this article we are interested in the Lie point symmetries, the Noether symmetries and the conservation laws of equation
\bb\label{1.3}
u_{tt}=u_{xx}+(\cot{x})\;u_{x}+\f{1}{\sin^{2}{x}}\;u_{yy}+f(u),
\ee
where $f:\R\rightarrow\R$ is a smooth function. Existence, uniqueness and behavior of solutions of initial value problems of (\ref{1.3}) are established in \cite{mo}.

Denoting by $Isom(\s)$ the Lie algebra of the Killing vector fields of $(\s,g)$, our main result can be formulated as follows:

\begin{theorem}\label{group}
Except to the linear cases, the symmetry Lie algebra of equation $(\ref{1.3})$ with an arbitrary function $f(u)$ is generated by $Isom(\s,g)$, that is,
\bb\label{s1}
S_{0}=\f{\p}{\p t},\,\,\,\,\,
S_{1}=\f{\p}{\p y},\,\,\,\,\,
S_{2}=\sin{y}\f{\p}{\p x}+\cot{x}\cos{y}\f{\p}{\p y},\,\,\,\,\,
S_{3}=\cos{y}\f{\p}{\p x}-\cot{x}\sin{y}\f{\p}{\p y}.
\ee

If $f(u)=c\,u$, $c=const.$, in $(\ref{1.2})$, in addition to $Isom(\s)$, we have the following generators:
\bb\label{s5}
S_{4}=u\f{\p}{\p u}
\ee
and
\bb\label{s}
S_{\infty}=b(x, y, t)\f{\p}{\p u},
\ee
where 
\bb\label{s6}
b_{tt}=b_{xx}+(\cot{x})\;b_{x}+\f{1}{\sin^{2}{x}}\;b_{yy}+c\,b.
\ee

The case $f(u)=k=const\neq 0$ is reduced to the case $f(u)=0$ under the change $u\rightarrow u+kt^{2}/2$.

\end{theorem}

As a consequence, we have the following classification of the Noether symmetries.

\begin{theorem}\label{noether}
For any function $f(u)$ in $(\ref{1.3})$, the isometry algebra $Isom(\s)$ generates a variational symmetry Lie algebra. If $f(u)=c\,u$, the symmetry $(\ref{s})$ is a Noether symmetry, where $b=b(x,y,t)$ satisfies $(\ref{s6})$.
\end{theorem}

From Theorem \ref{noether} and the Noether's Theorem, we have:

\begin{corollary}\label{conslaw}
Let $F=F(u)$ be a function such that $F'(u)=f(u)$. The conservation laws of the Noether symmetries of equation $(\ref{1.3})$, for any function $f(u)$, are:
\begin{enumerate}
\item For the symmetry $S_{0}$, the conservation law is $Div(A)=0$, where $A=(A^{0},A^{1},A^{2})$ and
\bb\label{cl1}
A^{0}  =  \ds{-\f{\sin{x}}{2}\,u_{t}^{2}-\f{\sin{x}}{2}\,u_{x}^{2}-\f{1}{2\sin{x}}\,u_{y}^{2}-\sin{x}F(u)},\,\,\,\,
A^{1}  =  \ds{\sin{x}\,u_{t}u_{x}},\,\,\,\,
A^{2}  =  \ds{\f{1}{\sin{x}}\,u_{t}u_{y}}.
\ee

\item For the symmetry $S_{1}$, the conservation law is $Div(B)=0$, where $B=(B^{0},B^{1},B^{2})$ and 
\bb\label{cl2}
B^{0}  =  \ds{-\sin{x}\,u_{t}u_{y}},\,\,
\,\,
B^{1}  =  \ds{\sin{x}\,u_{x}u_{y}},\,\,
\,\,
B^{2}  =  \ds{\f{\sin{x}}{2}\,u_{t}^{2}-\f{\sin{x}}{2}\,u_{x}^{2}+\f{1}{2\sin{x}}\,u_{y}^{2}-\sin{x}\,F(u)}.
\ee

\item For the symmetry $S_{2}$, the conservation law is $Div(C)=0$, where $C=(C^{0},C^{1},C^{2})$ and 
\bb\label{cl3}
\ba{l c l}
C^{0} & = & \ds{-\sin{x}\,\sin{y}\,u_{t}u_{x}-\cos{x}\,\cos{y}\,u_{t}u_{y}},\\
\\
C^{1} & = & \ds{\f{\sin{x}\,\sin{y}}{2}\,u_{t}^{2}+\f{\sin{x}\,\sin{y}}{2}\,u_{x}^{2}-\f{\sin{y}}{2\sin{x}}\,u_{y}^{2}+\cos{x}\,\cos{y}\,u_{x}\,u_{y}-\sin{x}\,\sin{y}\,F(u)},\\
\\
C_{2} & =& \ds{\f{\cos{x}\,\cos{y}}{2}\,u_{t}^{2}-\f{\cos{x}\,\cos{y}}{2}\,u_{x}^{2}+\f{\cos{x}\,\cos{y}}{2\sin^{2}{x}}\,u_{y}^{2}+\f{\sin{y}}{\sin{x}}\,u_{x}\,u_{y}-\cos{x}\,\cos{y}\,F(u)}.
\ea
\ee
\item For the symmetry $S_{3}$, the conservation law is $Div(D)=0$, where $D=(D^{0},D^{1},D^{2})$ and
\bb\label{cl4}
\ba{l c l}
D^{0} & = & \ds{-\sin{x}\,\cos{y}\,u_{t}\,u_{x}+\cos{x}\,\sin{y}\,u_{t}\,u_{y}},\\
\\
D^{1} & = &\ds{\f{\sin{x}\,\cos{y}}{2}\,u_{t}^{2}+\f{\sin{x}\,\cos{y}}{2}\,u_{x}^{2}-\f{\cos{y}}{2\sin{x}}\,u_{y}^{2}}\\
\\
&&\ds{-\cos{x}\,\sin{y}\,u_{x}\,u_{y}-\sin{x}\,\cos{y}\,F(u)},\\
\\ 
D^{2} & = & \ds{-\f{\cos{x}\,\sin{y}}{2}\,u_{t}^{2}+\f{\cos{x}\,\sin{y}}{2}\,u_{x}^{2}-\f{\cos{x}\,\sin{y}}{2\sin^{2}{x}}\,u_{y}^{2}}\\
\\
& &\ds{+\f{\cos{y}}{2\sin{x}}\,u_{x}\,u_{y}+\cos{x}\,\sin{x}\,F(u)}.
\ea
\ee
\item If $F(u)=c\,u^{2}/2$, then the conservation law for the symmetry $(\ref{s})$, with $b$ satisfying $(\ref{s6})$, is $Div(\al)=0$, where $\al=(\al^{0},\al^{1},\al^{2})$ and
\bb\label{cl5}
\al^{0}  =  \ds{\sin{x}\,(b\,u_{t}-b_{t}\,u)},\,\,
\,\,
\al^{1}  = \ds{\sin{x}\,(b_{x}\,u-b\,u_{x})},\,\,
\,\, 
\al^{2}  = \ds{\f{1}{\sin{x}}\,(b\,u_{t}-b_{t}\,u)}.
\ee

\end{enumerate}
\end{corollary}

We shall not present preliminaries concerning Lie point symmetries of differential equations supposing that the reader is familiar with the basic notions and methods of contemporary group analysis. See \cite{bk,ib,ol}. For a geometric viewpoint of Lie point symmetries, see \cite{man,olvjdg}.

This paper is organized as follows. In section \ref{s2} we recall some geometric results regarding to $(\s,g)$. These results will be used in section \ref{gc} to prove the Theorem 1. The Noether's symmetries and the conservation laws are obtained in section \ref{no}. In section \ref{al} we identify the classical Lie algebras that the symmetry Lie algebras are isomorphic to. 

\section{The product manifold $\s$}\label{s2}

\

Let $x^{0}=t,\,x^{1}=x$ and $x^{2}=y$ be local coordinates of $(\s,g)$. We observe that the Riemann and the Ricci tensors used in this paper coincide with those in Yano's book \cite{ya} and in Dubrovin, Fomenko and Novikov's book \cite{dnf}.

\begin{lemma}\label{curves}
The scalar curvature $R$ of the product manifold $(\s,g)$ is constant. 
\end{lemma}

\begin{proof}
The Riemann tensor of $(\s,g)$ is $$R^{i}_{jks}=-(\de^{2i}\de_{1j}\de_{2k}\de_{1s}-\de^{2i}\de_{1j}\de_{1k}\de_{2s})+\sin^{2}{x}(\de^{2i}\de_{2j}\de_{2k}\de_{1s}-\de^{1i}\de_{2j}\de_{1k}\de_{2s}).$$ Then $R^{i}_{s}=-\de^{2i}\de_{2s}-\de^{2i}\de_{1s}$ and $R=-1$.
\end{proof}

\begin{lemma}\label{curvsec}
The sectional curvature of $(\s,g)$ is non-constant.
\end{lemma}

\begin{proof}
Let $K(p,X,Y)$ be the sectional curvature of $(\s,g)$ at a point $p=(t,x,y)$. (See \cite{lov} for the definitions.) Let $X=(X^{0},X^1,X^2)$ and $Y=(Y^0,Y^1,Y^2)$. Then, we obtain
$$K(p,X,Y)=\f{-X^{2}\,\sin^{2}{x}+X^1}{2X^1}.$$
\end{proof}

\begin{lemma}\label{iso}
The isometry group of $(\s,g)$ is generated by the vector fields $S_{0},\,S_{1},\,S_{2}$ and $S_{3}$.

\end{lemma}

\begin{proof}
It is clear that the vector fields (\ref{s1}) satisfy the equation 
$$L_{X}g_{ij}=\xi^{s}\f{\p g_{ij}}{\p x^{s}}+g_{kj}\f{\p \xi^{k}}{\p x^{i}}+g_{ik}\f{\p \xi^{k}}{\p x^{j}}=0.$$

From Lemma \ref{curvsec}, the sectional curvature of $(\s,g)$ is non-constant. Then, from Yano \cite{ya}, pag. 57, Theorem 6.2, $dim(Isom(\s))< 6$. From Fubini's Theorem (see Yano \cite{ya}), $dim(Isom(\s))$ cannot be $5$. Thus, $dim(Isom(\s))\leq 4$. Since (\ref{s1}) are isometries, we conclude that the isometry algebra $Isom(\s)$ is generated by $S_{0},\,S_{1},\,S_{2}$ and $S_{3}$.
\end{proof}

\section{The group classification}\label{gc}

\

In this section we perform the group classification of equation (\ref{1.3}). To begin with we need of the following Lemma 

\begin{lemma}\label{igor}
Let $(M^{n},g)$ be a manifold with non-null constant scalar curvature. Then the Lie point symmetry group of the Poisson equation $(\ref{1.2})$ with an arbitrary $f(u)$ coincides with the isometry group of $M^{n}$. 

In the particular cases $f(u)=c\,u$, where $c=const.$, in addition to the isometry group, we have the generators 
$$
\ba{l c l}
\ds{U=u\f{\p}{\p u} }& \text{ and } & \ds{X_{\infty}=b(x)\f{\p}{\p u}},
\ea
$$
where $b$ satisfies $(\ref{1.2})$.
\end{lemma}
\begin{proof} See \cite{yi1}.\end{proof}
We observe that the Lemma \ref{igor} is a particular case of the main result obtained in \cite{yi1}. In this work the authors carried out the group classification of equation (\ref{1.2}) on an arbitrary (pseudo) Riemannian manifold. 

{\bf Proof of Theorem 1}: From Lemma \ref{curves}, the scalar curvature of $(\s,g)$ is $R=-1$. Then, the Theorem \ref{group} follows from lemmas \ref{iso} and \ref{igor}.\qed

\section{The Noether's symmetries and the conservation laws}\label{no}

\

In this section we prove Theorem \ref{noether}. 

It is easy to check that the if $X\in Isom(\s)$, then $X$ is a variational symmetry of equation (\ref{1.2}), for any function $f(u)$. That is $$X^{(1)}{\cal L}+{\cal L} D_{i}\xi^{i}=0,$$ where
\bb\label{3.1}
{\cal L}=\f{\sin{x}}{2}u_{t}^{2}-\f{\sin{x}}{2}u_{x}^{2}-\f{1}{2\sin{x}}u_{y}^{2}+\sin{x}\,F(u)
\ee
is the function of Lagrange of equation (\ref{1.2}). For more details, see \cite{yi1}.

Let us consider the symmetry (\ref{s}). It is easy to verify that 
$$X^{(1)}{\cal L}+{\cal L} D_{i}\xi^{i}=Div\left(\sin{x}\,b_{t}\,u,-\sin{x}\,b_{x}\,u,-\f{1}{\sin{x}}\,b_{y}\,u\right),$$
where $F(u)=c\,u^{2}/2$ in (\ref{3.1}). 
These observations prove the Theorem \ref{noether}.

The following lemma establishes the conservation laws (\ref{cl1})-(\ref{cl5}):

\begin{lemma}\label{consl} The conservation laws of the Noether symmetries of equation $(\ref{1.2})$, where $(M^{n},g)$ is a manifold with constant, non-null scalar curvature, are $D_{i}A^{i}=0$, where
\bb\label{cc1}
A^{k}=\sqrt{g}(\f{1}{2}g^{ij}\xi^{k}-g^{kj}\xi^{i})u_{i}u_{j}-\sqrt{g}\xi^{k}F(u),\ee
for any function $f(u)$. If $f(u)=c\,u$, then the conservation law corresponding to the symmetry $(\ref{s})$ is
\bb\label{cc5}
A^{k}   = \sqrt{g}g^{jk}(bu_{j}-b_{j}u).
\ee
\end{lemma}
\begin{proof} It is a consequence of \cite{yi1} when the scalar curvature of $(M^{n},g)$ is constant.\end{proof}
{\bf Proof of Corollary \ref{conslaw}}: Substituting the symmetries and the metric coefficients into (\ref{cc1})-(\ref{cc5}), we obtain (\ref{cl1})-(\ref{cl5}).\qed

\section{Symmetry Lie algebras}\label{al}

\

Let $\mathfrak{S}_{1},\,\mathfrak{S}_{2}$ be the finite dimensional symmetry Lie algebras for an arbitrary $f(u)$ and $f(u)=c\,u$, $c=const$, respectively. Following the notations of \cite{pswz,pw}, the symmetry Lie algebras are:

\begin{enumerate}
\item If $f(u)$ is an arbitrary function, then $[S_{1},S_{2}]=S_{3},\,[S_{1},S_{3}]=-S_{2},\,[S_{2},S_{3}]=S_{1}$. Thus, $\mathfrak{S}_{1}=Isom(\s)\approx A_{3,9}\oplus A_{1}$, where $A_{3,9}=\mathfrak{s0}(3)$. 
\item If $f(u)=c\,u$, then $\mathfrak{S}_{2}\approx A_{3,9}\oplus 2A_{1}$.
\end{enumerate}

We have the following one-dimensional subálgebras of $\mathfrak{S}_{1}$: $\mathfrak{L}_{1}=<S_{0}+a\,S_{1}>$ and $\mathfrak{L}_{2}=<S_{1}>$.

If $f(u)=c\,u$, we have the following classification of subalgebras of $\mathfrak{S}_{2}$:
\begin{enumerate}
\item Dimension 1: $\mathfrak{L}_{1}=<a\,S_{0}+S_{1}+b\,S_{4}>\approx A_{1}$ and $\mathfrak{L}_{2}=<a\,S_{0}+b\,S_{4}>\approx A_{1}$.
\item Dimension 2: $\mathfrak{L}_{3}=<a\,S_{0}+b\,S_{4},S_{1}>\approx 2A_{1}$ and $\mathfrak{L}_{4}=<S_{0},S_{4}>\approx 2A_{1}$.
\item Dimension 3: $\mathfrak{L}_{5}=<S_{1},S_{2},S_{3}>\approx A_{3,9}$ and $\mathfrak{L}_{6}=<S_{0},S_{1},S_{4}>\approx 3A_{1}$.
\item Dimension 4: $\mathfrak{L}_{7}=<a\,S_{0}+b\,S_{4}, S_{1},S_{2},S_{3}>\approx A_{3,9}\oplus A_{1}$.
\end{enumerate}

We observe that the subalgebras (1-4) above were obtained by Azad and Mustafa when $f(u)=0$ in (\ref{1.3}). 

The invariant solutions of (\ref{1.3}) can be obtained following the same procedure employed by Azad and Mustafa in \cite{A} with addition of the corresponding term $f(u)$ in (\ref{1.2}). Thus we shall omit the details. \\
\\

{\bf Acknowledgements}\\
\\

I am grateful to Y. Bozhkov for his careful reading of this paper as well as for his firm encouragement. I am also pleased to thank the anonymous referee for his comments.

\end{document}